\documentclass[11pt, reqno, draft]{amsart}

\usepackage{amsmath}
\usepackage{amssymb}
\usepackage{amsthm}
\usepackage{mathrsfs}
\usepackage{mathtools}
\usepackage{graphics}
\usepackage{bm}
\usepackage{comment}
\usepackage{color}

\usepackage[dvipdfmx, hypertexnames=false, hidelinks]{hyperref}
\usepackage{cleveref}

\usepackage{autonum}

\hypersetup{final}

\allowdisplaybreaks[3]
\theoremstyle{plain} 
\newtheorem{theorem}{Theorem}
\newtheorem{lemma}{Lemma}[section]

\newtheorem*{conjecture*}{Conjecture}
\newtheorem*{theorem*}{Theorem}
\newtheorem*{question*}{Question}

\theoremstyle{plain}

\theoremstyle{remark}
\newtheorem{remark}{Remark}

\theoremstyle{definition}
\newtheorem*{assumption*}{Assumption}

\newtheorem*{notations*}{Notations}
\newtheorem*{acknowledgment*}{Acknowledgments}

\numberwithin{equation}{section}

\crefname{section}{Section}{Sections}

\crefname{theorem}{Theorem}{Theorems}
\crefname{corollary}{Corollary}{Corollaries}
\crefname{lemma}{Lemma}{Lemmas}
\crefname{proposition}{Proposition}{Propositions}
\crefname{claim}{Claim}{Claims}
\crefname{definition}{Definition}{Definitions}
\crefname{notation}{Notation}{Notations}
\crefname{problem}{Problem}{Problems}
\crefname{note}{Note}{Notes}
\crefname{remark}{Remark}{Remarks}
\crefname{example}{Example}{Examples}
\crefname{equation}{}{}

\crefname{enumi}{}{}
\crefname{enumii}{}{}
\crefname{enumiii}{}{}

\newcommand\swapcommand[2]{%
\let\swaptemp#1
\let#1#2
\let#2\swaptemp
}

\let\sl\l
\renewcommand\l{%
		\leavevmode
	\ifmmode
	\left
	\else
		\sl
	\fi
}

\let\sL\L
\renewcommand\L{%
		\leavevmode
	\ifmmode
	\mathscr{L}
	\else
		\sL
	\fi
}

\swapcommand{\SS}{\S}
\renewcommand{\S}{\mathscr{S}}

\newcommand{\ZZ}{\mathbb{Z}}

\newcommand{\lam}{\lambda}

\newcommand{\ceq}{\coloneqq}

\renewcommand{\a}{\alpha}
\renewcommand{\b}{\beta}
\renewcommand{\i}{\mathrm{i}}

\renewcommand{\r}{\right}

\renewcommand{\Re}{\operatorname{Re}}
\renewcommand{\Im}{\operatorname{Im}}
\renewcommand{\epsilon}{\varepsilon}

\title{A note on $r$-gaps between zeros of the Riemann zeta-function}
\author[S. Inoue]{Sh\={o}ta Inoue}

\address{Department of Mathematics, Tokyo Institute of Technology, 2-12-1 Ookayama, Meguro-ku, Tokyo 152-8551, Japan}
\email{inoue.s.bd@m.titech.ac.jp}

\keywords{The Riemann zeta-function, Gaps of zeros of the Riemann zeta-function, The Riemann Hypothesis}

\subjclass{Primary 11M06; Secondary 11M26}

\begin{document}

\begin{abstract}
	In this note, we prove Selberg's announced result on $r$-gaps between zeros of the Riemann zeta-function $\zeta$.
	Our proof uses a result on variations of $\arg\zeta$ by Tsang based on Selberg's method.
	The same result with explicit constants under the Riemann Hypothesis has been obtained by Conrey and Turnage-Butterbaugh using a different method.
	We explain how to obtain explicit constants under the Riemann Hypothesis using our approach which is based on Selberg's and Tsang's arguments.

\end{abstract}

\maketitle

\section{\textbf{Introduction}}

	Let $N(T)$ be the number of zeros $\rho = \b + \i \gamma$ of the Riemann zeta-function $\zeta$ with $0 < \gamma \leq T$ counted with multiplicity,
	where if $T$ is equal to the imaginary part of some zeros, then the multiplicity is to be counted with weight $1/2$.
	Denote the sequence of ordinates of zeros of $\zeta$ on the upper half plane by
	$0 < \gamma_{1} \leq \gamma_{2} \leq \cdots \leq \gamma_{n} \leq \gamma_{n + 1} \leq \cdots$.
	The Riemann-von Mangoldt formula
	\begin{align}
		\label{RVFEx}
		N(T)
		&= \frac{T}{2\pi}\log\l( \frac{T}{2\pi e} \r)
		+ S(T) + \frac{7}{8} + O\l(\frac{1}{T}\r)\\
		\label{RVFEr}
		&\sim \frac{T}{2\pi}\log{T}
	\end{align}
	is known as a useful formula for understanding the distribution of zeros.
	Here, $S(t) \ceq \frac{1}{\pi}\Im \log{\zeta(\frac{1}{2} + \i t)} = \frac{1}{\pi}\Im \int_{\infty}^{1/2}\frac{\zeta'}{\zeta}(\a + \i t)d\a$ if $t$ is not equal to the imaginary part of zeros of $\zeta$.
	If $t$ is equal to the imaginary part of some zero, then $S(t) = (S(t + 0) + S(t - 0)) / 2$.
	Define
	\begin{align}
		\lam_{r} \ceq \limsup_{n \rightarrow +\infty}\frac{\gamma_{n + r} - \gamma_{n}}{2\pi r / \log{\gamma_{n}}}
		\quad \text{ and } \quad
		\mu_{r} \ceq \liminf_{n \rightarrow +\infty}\frac{\gamma_{n + r} - \gamma_{n}}{2\pi r / \log{\gamma_{n}}}.
	\end{align}
	We have the trivial bound $\mu_{r} \leq 1 \leq \lam_{r}$ coming from \cref{RVFEr}.
	In the case $r = 1$, there are many results known concerning the nontrivial bound $\mu_{1} < 1 < \lam_{1}$.
	Montgomery \cite{Mo1973} showed such a nontrivial bound under the Riemann Hypothesis.
	In his work, he suggested a conjecture called the pair correlation conjecture.
	The pair correlation conjecture suggests that $\lam_{1} = +\infty$ and $\mu_{1} = 0$.
	Explicit bounds of $\lam_{1}$ and $\mu_{1}$ are known such as \cite{BM2018}, \cite{BMN2010}, \cite{CGG1984}, \cite{Ha2005}, \cite{Mu1982}, \cite{Pr2016}.
	The best bounds known on the Riemann Hypothesis are $\lam_{1} > 3.18$ and $\mu_{1} < 0.515396$ given by Bui and Milinovich \cite{BM2018} and Preobrazhenski\u{i} \cite{Pr2016} respectively.
	Prior to these, Selberg \cite{Se1946zRH} announced without any proof that the nontrivial bound $\mu_{1} < 1 < \lam_{1}$ holds unconditionally, which was again mentioned by Fujii \cite{Fu1975}.
	Later, Heath-Brown \cite[Section 9.26]{TRZ} gave a proof of their announced result by using the second moment of $S(t)$ in short intervals due to Fujii \cite{Fu1974}.

	In this note, we consider the bounds of $\lam_{r}$ and $\mu_{r}$ for general $r$.
	Selberg \cite{Se1946zRH}, \cite{SCP} and Fujii \cite{Fu1975} discussed this topic, and
	Fujii stated that $\lam_{r} > 1 + \exp(-C r^{2})$ and $\mu_{r} < 1 - \exp(-C r^{2})$ hold for each $r \in \ZZ_{\geq 1}$ with $C$ a positive absolute constant.
	Also, Selberg \cite[p.\ 355]{SCP} announced that $\lam_{r}$ and $\mu_{r}$ satisfy the sharper inequalities
	\begin{align}
		\label{SeARRG}
		\lam_{r} > 1 + \frac{\theta}{r^{\a}}
		\quad \text{ and } \quad
		\mu_{r} < 1 - \frac{\theta}{r^{\a}}
	\end{align}
	for all positive integer $r$ with $\theta$ an absolute positive constant,
	where $\a$ may be taken as $2 / 3$ unconditionally, and if we assume the Riemann Hypothesis as $1 / 2$.
	However, they did not give any proof for those inequalities.

	Recently, Conrey and Turnage-Butterbaugh \cite{CTB2018} gave a proof of the conditional result announced by Selberg.
	However, the unconditional statement of Selberg has not yet been proven.
	One of the purposes of this note is to give a proof of it.
	Actually, we prove the following theorem.

	\begin{theorem}	\label{UR_rgaps}
		There exist absolute positive constants $\theta_{1}, \theta_{2}$ such that
		for any $r \in \ZZ_{\geq 1}$
		\begin{align}	\label{INE_UR_rgaps}
			\limsup_{n \rightarrow +\infty}\frac{\gamma_{n + r} - \gamma_{n}}{2\pi r / \log{\gamma_{n}}} > 1 + \frac{\theta_{1}}{r^{2/3}}
			\quad \text{ and } \quad
			\liminf_{n \rightarrow +\infty}\frac{\gamma_{n + r} - \gamma_{n}}{2\pi r / \log{\gamma_{n}}} < 1 - \frac{\theta_{2}}{r^{2/3}}.
		\end{align}
	\end{theorem}

	In the next section, we prove this theorem.
	Later in this note, we also discuss the conditional bounds of $\lam_{r}$ and $\mu_{r}$ and improve the result by Conrey and Turnage-Butterbaugh.

\section{\textbf{Proof of \cref{UR_rgaps}}}

	In the proof of \cref{UR_rgaps}, we only use the Riemann-von Mangoldt formula \cref{RVFEx}
	and the $\Omega$-result of $S(T)$ in short intervals shown by Tsang \cite{Ts1986}.

	\begin{lemma}	\label{SOmega}
		There exists an absolute positive constant $c$ such that for any large $T$ and $h \in [\frac{1}{\log{T}}, \frac{1}{\log{\log{T}}}]$
		\begin{align}
			\sup_{t \in [T, 2T]}\l\{\pm(S(t + h) - S(t))\r\}
			\geq c (h \log{T})^{1/3}.
		\end{align}
	\end{lemma}

	Assuming the Riemann Hypothesis, we can improve the exponent $1/3$ in \cref{SOmega} to $1/2$.
	Using it, we can also prove Selberg's announced conditional result.

	\begin{proof}[Proof of \cref{UR_rgaps}]
		Since we can prove the first and second assertion by the same argument, we only give the proof of the first assertion.
		Let $T$ be a sufficiently large number, and $h$ be a small number to be chosen later.
		If there exists some $t \in [T, 2T]$ such that
		\begin{align}	\label{EIE}
			N(t + h) - N(t) < r,
		\end{align}
		then we have
		\begin{align}
			\sup_{\gamma_{n}, \gamma_{n + r} \in [T, 2T + h]}\frac{\gamma_{n + r} - \gamma_{n}}{h} > 1.
		\end{align}
		By the Riemann-von Mangoldt formula, we have
		\begin{align}
			N(t + h) - N(t)
			= \frac{h}{2\pi}\log{t} + S(t + h) - S(t) + o(1)
		\end{align}
		for any $t \in [T, 2T]$ as $h \rightarrow 0$ and $T \rightarrow +\infty$.
		We use this and \cref{SOmega} to obtain
		\begin{align}
			N(t + h) - N(t)
			\leq \frac{h}{2\pi}\log{t} - c(h \log{t})^{1/3} + o(1)
		\end{align}
		for some $t \in [T, 2T]$. For this $t$, if $h$ satisfies the equation
		\begin{align}
			\frac{h}{2\pi r / \log{T}} = 1 + \frac{\theta_{1}}{r^{2/3}}
		\end{align}
		with $\theta_{1}$ a suitable positive constant, then inequality \cref{EIE} holds.
		Hence, we find that
		\begin{align}
			&\limsup_{n \rightarrow + \infty}\frac{\gamma_{n + r} - \gamma_{n}}{2\pi r / \log{\gamma_{n}}}
			\geq \lim_{T \rightarrow + \infty}\sup_{\gamma_{n}, \gamma_{n + r} \in [T, 2T + h]}
			\frac{\gamma_{n + r} - \gamma_{n}}{2\pi r / \log{\gamma_{n}}}\\
			&= \lim_{T \rightarrow + \infty}\sup_{\gamma_{n}, \gamma_{n + r} \in [T, 2T + h]}
			\frac{\gamma_{n + r} - \gamma_{n}}{h}\frac{h}{2\pi r / \log{T}}\frac{\log{\gamma_{n}}}{\log{T}}
			> 1 + \frac{\theta_{1}}{r^{2/3}},
		\end{align}
		which completes the proof of \cref{UR_rgaps}.
	\end{proof}

\section{\textbf{Conditional result}}	\label{Sec_CR}

	Conrey and Turnage-Butterbaugh \cite{CTB2018} gave a proof of Selberg's announced conditional result.
	They also discussed the size of $\theta$ in \cref{SeARRG}.
	In fact, they proved under the Riemann Hypothesis that
	\begin{align}	\label{SARRH}
		\limsup_{n \rightarrow +\infty}\frac{\gamma_{n + r} - \gamma_{n}}{2\pi r / \log{\gamma_{n}}} > 1 + \frac{\Theta}{r^{1/2}}
		\quad \text{ and } \quad
		\liminf_{n \rightarrow +\infty}\frac{\gamma_{n + r} - \gamma_{n}}{2\pi r / \log{\gamma_{n}}} < 1 - \frac{\vartheta}{r^{1/2}}
	\end{align}
	for any positive integer $r$ with $\Theta = 0.599648$ and $\vartheta = 0.379674$
	(see Remark in \cite{CTB2018}).
	Moreover, they also proved that inequalities \cref{SARRH} hold with $\Theta, \vartheta = A_{0} + o(1)$ as $r \rightarrow + \infty$.
	Here, $A_{0} \ceq \max_{B > 0}\frac{2B}{\pi}\arctan\l( \frac{\pi}{B^{2}} \r) = 0.9064997 \cdots$.

	Their method is based on the works of Montgomery and Odlyzko \cite{MO1984} and Conrey, Ghosh, and Gonek \cite{CGG1984}.
	Here, we apply the method of Selberg and Tsang used in the proof of \cref{UR_rgaps},
	and compare those two distinct methods.
	Using the argument of Tsang in \cite{Ts1986}, we can prove \cref{SARRH} with $\Theta = 0.414269$ and $\vartheta = 0.403816$ for any $r \geq 5$ under the Riemann Hypothesis.
	Hence, we can improve partially the result of Conrey and Turnage-Butterbaugh
	for $\vartheta$ by the method of Tsang.
	A sketch of this fact will be described in the final section of the present paper.
	On the other hand, we can obtain a further improvement by adapting the method of Conrey and Turnage-Butterbaugh and by changing their way of calculating certain integrals.
	Roughly speaking, they calculated the integrals defined by \cref{CTBSCT,CTBSCvt} below by dividing into smaller integrals.
	On the other hand, we evaluate the integrals by using the saddle point method.
	As a result, we show the following theorem.

	\begin{theorem}	\label{CR_rgaps}
		Under the Riemann Hypothesis, inequalities \cref{SARRH} hold with the constants $\Theta = A_{0} = 0.9064997 \cdots$ and $\vartheta = 0.484604$ uniformly for $r \geq 1$.
	\end{theorem}

	\begin{remark}
		The constant of $\vartheta$ in \cref{CR_rgaps} comes from the case $r = 1$.
		In the case $r = 1$, we use the result of Preobrazhenski\u{i} \cite{Pr2016}.
		For $r \geq 2$, we prove \cref{SARRH} with $\vartheta = 0.61861$.
	\end{remark}

	\begin{proof}
		When $r = 1$, this theorem has been proven by Bui and Milinovich \cite{BM2018} for $\Theta = 2.18$, and by Preobrazhenski\u{i} \cite{Pr2016} for $\vartheta = 0.484604$ respectively.
		Therefore, we consider only the case $r \geq 2$.
		We use the same sufficient condition for \cref{SARRH} as in the work of Conrey and Turnage-Butterbaugh \cite{CTB2018}.
		In fact, by equations (2.1) and (2.7) in \cite{CTB2018} and the explanations written below the equations, the statement that
		there exists a real number $\ell \geq 1$ such that
		\begin{align}
			\label{CTBSCT}
			\Theta
			&< \frac{2\ell}{\sqrt{r}}\int_{0}^{1}\frac{\sin(\pi(r + \Theta \sqrt{r})v)}{\pi v}(1 - v)^{\ell^{2}}dv,\\
			\label{CTBSCvt}
			\vartheta
			&< \frac{2\ell}{\sqrt{r}}\int_{0}^{1}\frac{\sin(\pi(r - \vartheta \sqrt{r}) v)}{\pi v}(1 - v)^{\ell^{2}}dv
		\end{align}
		is a sufficient condition of \cref{SARRH} for every $r$.
		Note that we have taken the limits $T \rightarrow + \infty$, $\delta \rightarrow 0$ for the original argument in \cite{CTB2018}.
		We can confirm by numerical calculations that \Cref{TBTvtCT} presents approximations of $\ell$, $\Theta$, and $\vartheta$ satisfying \cref{CTBSCT,CTBSCvt} for $1 \leq r \leq 20$.
		\begin{table}[htb]
			\centering
			\caption{Approximations of $(\Theta, \ell)$ and $(\vartheta, \ell)$ satisfying \cref{CTBSCT,CTBSCvt}, which imply \cref{SARRH}}
			\label{TBTvtCT}
			\small
			\scalebox{1}[1]{
				\begin{tabular*}{140mm}{@{\extracolsep{\fill}}lcr|lcr}
					\hline
					$r$	&	$(\Theta, \ell)$	&	$(\vartheta, \ell)$	&	$r$	&	$(\Theta, \ell)$	&	$(\vartheta, \ell)$\\
					\hline
					1	&	(1.337, 2.16)	&	(0.482, 1.02)			&	11	&	(1.032, 5.66)	&	(0.784, 4.29)\\
					\hline
					2	&	(1.208, 2.8)	&	(0.6186133, 1.41963)	&	12	&	(1.027, 5.88)	&	(0.789, 4.51)\\
					\hline
					3	&	(1.151, 3.28)	&	(0.675, 1.88)			&	13	&	(1.022, 6.09)	&	(0.794, 4.73)\\
					\hline
					4	&	(1.117, 3.68)	&	(0.706, 2.29)			&	14	&	(1.018, 6.3)	&	(0.798, 4.93)\\
					\hline
					5	&	(1.094, 4.03)	&	(0.727, 2.66)			&	15	&	(1.014, 6.49)	&	(0.802, 5.13)\\
					\hline
					6	&	(1.078, 4.36)	&	(0.742, 2.98)			&	16	&	(1.0107, 6.69)	&	(0.805, 5.32)\\
					\hline
					7	&	(1.065, 4.65)	&	(0.754, 3.28)			&	17	&	(1.007, 6.87)	&	(0.808, 5.51)\\
					\hline
					8	&	(1.054, 4.93)	&	(0.764, 3.56)			&	18	&	(1.004, 7.05)	&	(0.811, 5.69)\\
					\hline
					9	&	(1.046, 5.18)	&	(0.772, 3.81)			&	19	&	(1.002, 7.22)	&	(0.813, 5.86)\\
					\hline
					10	&	(1.038, 5.43)	&	(0.778, 4.06)			&	20	&	(0.9995, 7.39)	&	(0.815, 6.03)\\
					\hline
				\end{tabular*}
			}
		\end{table}
		Hence, \cref{CR_rgaps} holds for $1 \leq r \leq 7$.

		In the following, we assume that $r \geq 8$ and $\ell \geq 2$.
		For $\theta \in [-1, 1]$, let
		\begin{align}
			F(r)
			\ceq \frac{2\ell}{\sqrt{r}}\int_{0}^{1}\frac{\sin(\pi(r + \theta \sqrt{r}) v)}{\pi v}(1 - v)^{\ell^{2}}dv.
		\end{align}
		Since the inequality $(1 - v)^{\ell^{2}} \leq e^{-\ell^2 v}$ holds for $0 \leq v \leq 1$, we have
		\begin{align}
			\bigg|\int_{\ell^{-1}}^{1}\frac{\sin(\pi(r + \theta \sqrt{r}) v)}{\pi v}(1 - v)^{\ell^{2}}dv\bigg|
			\leq \int_{\ell^{-1}}^{1}\frac{\ell}{\pi}e^{-\ell^2 v}dv
			\leq \frac{1}{\pi\ell}e^{-\ell}.
		\end{align}
		We write
		\begin{align}
			&\int_{0}^{\ell^{-1}}\frac{\sin(\pi(r + \theta \sqrt{r}) v)}{\pi v}(1 - v)^{\ell^{2}}dv\\
			&= \int_{0}^{\ell^{-1}}\frac{\sin(\pi(r + \theta \sqrt{r}) v)}{\pi v}e^{-\ell^{2} v}dv
			- \int_{0}^{\ell^{-1}}\frac{\sin(\pi(r + \theta \sqrt{r}) v)}{\pi v}\l\{ e^{-\ell^{2} v} - (1 - v)^{\ell^{2}} \r\}dv.
		\end{align}
		Since the inequalities
		\begin{align}
			0 \leq e^{-\ell^{2} v} - (1 - v)^{\ell^{2}}
			= e^{-\ell^{2}v}\l( 1 - \exp\l( -\ell^{2}\sum_{n = 2}^{\infty}\frac{v^{n}}{n} \r) \r)
			\leq \ell^{2}e^{-\ell^{2} v}\sum_{n = 2}^{\infty}\frac{v^{n}}{n}
			\leq \ell^{2}e^{-\ell^{2}v} v^{2}
		\end{align}
		hold for $0 \leq v \leq \frac{1}{\ell} \leq \frac{1}{2}$, we obtain
		\begin{align}
			\bigg|\int_{0}^{\ell^{-1}}\frac{\sin(\pi(r + \theta \sqrt{r}) v)}{\pi v}\l\{ e^{-\ell^{2} v} - (1 - v)^{\ell^{2}} \r\}dv \bigg|
			&\leq \ell^{2}(r + \theta \sqrt{r})\int_{0}^{\infty} e^{-\ell^{2}v} v^{2} dv\\
			&= \frac{2(r + \theta \sqrt{r})}{\ell^{4}}.
		\end{align}
		We also write
		\begin{align}
			&\int_{0}^{\ell^{-1}}\frac{\sin(\pi(r + \theta\sqrt{r})v)}{\pi v}e^{-\ell^{2} v}dv\\
			&= \int_{0}^{\infty}
			\frac{\sin(\pi(r + \theta\sqrt{r})v)}{\pi v}e^{-\ell^{2} v}dv
			- \int_{\ell^{-1}}^{\infty}
			\frac{\sin(\pi(r + \theta\sqrt{r})v)}{\pi v}e^{-\ell^{2} v}dv.
		\end{align}
		The former integral is equal to $\frac{\arctan(\pi (r + \theta \sqrt{r}) / \ell^{2})}{\pi}$, which can be obtained by considering the Taylor expansions of $\sin$ and of $\arctan$.
		Also, the absolute value of the latter integral is
		$
			\leq \frac{1}{\pi\ell}e^{-\ell}.
		$
		From the above inequalities, we have
		\begin{align}
			\label{INE_Fr_1}
			&F(r) \geq
			\frac{2\ell}{\sqrt{r}}\l( \frac{1}{\pi}\arctan\l( \frac{\pi (r + \theta \sqrt{r})}{\ell^{2}}\r)
			- \frac{2e^{-\ell}}{\pi \ell} - \frac{2(r + \theta \sqrt{r})}{\ell^{4}} \r).
		\end{align}
		Here, we choose $\ell = B_{0}\sqrt{r + \theta \sqrt{r}}$ with $B_{0} = 1.502432\cdots$ the maximum point of the function $B \arctan(\pi / B^{2})$ in $B > 0$.
		Note that our assumption $\ell \geq 2$ holds for this $\ell$ when $r \geq 8$ and $|\theta| \leq 1$.
		Then inequality \cref{INE_Fr_1} is rewritten to
		\begin{align}
			F(r)
			&\geq
			2B_{0}\sqrt{1 + \tfrac{\theta}{\sqrt{r}}}\l( \frac{1}{\pi}\arctan\l( \frac{\pi}{B_{0}^{2}} \r)
			- \frac{2e^{-B_{0}\sqrt{r + \theta \sqrt{r}}}}{\pi B_{0}\sqrt{r + \theta \sqrt{r}}} - \frac{2}{B_{0}^{4} (r + \theta \sqrt{r})} \r)\\
			&= A_{0}
			+ \Biggl\{\l(\sqrt{1 + \tfrac{\theta}{\sqrt{r}}} - 1\r)A_{0}
			- \frac{4}{\sqrt{r}}\l( \frac{e^{-B_{0}\sqrt{r + \theta \sqrt{r}}}}{\pi} + \frac{1}{B_{0}^{3} \sqrt{r + \theta \sqrt{r}}} \r) \Biggr\}\\
			&= A_{0}
			+ \frac{1}{\sqrt{r}}\Biggl\{\frac{\theta}{\sqrt{1 + \frac{\theta}{\sqrt{r}}} + 1}A_{0} - 4\l( \frac{e^{-B_{0}\sqrt{r + \theta \sqrt{r}}}}{\pi } + \frac{1}{B_{0}^{3} \sqrt{r + \theta \sqrt{r}}} \r) \Biggr\}.
		\end{align}
		Hence, we see using the third and second equations that the inequalities
		\begin{align}
			\label{FIET}
			\Theta
			&< A_{0}
			+ \frac{1}{\sqrt{r}}\Biggl\{\frac{\Theta}{\sqrt{1 + \frac{\Theta}{\sqrt{r}}} + 1}A_{0} - 4\l( \frac{e^{-B_{0}\sqrt{r + \Theta \sqrt{r}}}}{\pi } + \frac{1}{B_{0}^{3} \sqrt{r + \Theta \sqrt{r}}} \r) \Biggr\},\\
			\label{FIEvt}
			\vartheta
			&< A_{0}
			- \Biggl\{\l(1 - \sqrt{1 - \tfrac{\vartheta}{\sqrt{r}}}\r)A_{0}
			+ \frac{4}{\sqrt{r}}\l( \frac{e^{-B_{0}\sqrt{r - \vartheta \sqrt{r}}}}{\pi} + \frac{1}{B_{0}^{3} \sqrt{r - \vartheta \sqrt{r}}} \r) \Biggr\}
		\end{align}
		imply inequalities \cref{CTBSCT}, \cref{CTBSCvt} respectively.

		For $\Theta = A_{0}$, the content of the braces $\{\}$ on the right hand side of \cref{FIET} is monotonically increasing for $r \geq 1$, and it holds that
		\begin{align}
			\frac{\Theta}{\sqrt{1 + \frac{\Theta}{\sqrt{8}}} + 1}A_{0}
			- 4\l( \frac{e^{-B_{0}\sqrt{8 + \Theta \sqrt{8}}}}{\pi} + \frac{1}{B_{0}^{3} \sqrt{8 + \Theta \sqrt{8}}} \r)
			\geq 0.009 > 0.
		\end{align}
		Therefore, when $\Theta = A_{0}$, inequality \cref{FIET} holds uniformly for $r \geq 8$.
		Combining this and the case $1 \leq r \leq 7$ we obtain the assertion of \cref{CR_rgaps} for $\Theta$.

		For $\vartheta = 0.61861$, the content of the braces $\{\}$ on the right hand side of \cref{FIEvt} is monotonically decreasing for $r \geq 1$, and it holds that
		\begin{align}
			A_{0}
			- \Biggl\{\l(1 - \sqrt{1 - \tfrac{\vartheta}{\sqrt{8}}}\r)A_{0}
			+ \frac{4}{\sqrt{8}}\l( \frac{e^{-B_{0}\sqrt{8 - \vartheta \sqrt{8}}}}{\pi} + \frac{1}{B_{0}^{3} \sqrt{8 - \vartheta \sqrt{8}}} \r) \Biggr\}
			\geq 0.62 > 0.61861.
		\end{align}
		Hence, when $\vartheta = 0.61861$, inequality \cref{FIEvt} holds uniformly for $r \geq 8$.
		By this and the results for $1 \leq r \leq 7$, we also obtain the assertion of \cref{CR_rgaps} for $\vartheta$.
	\end{proof}

\section{\textbf{Concluding remark}}
	As mentioned at the beginning of \cref{Sec_CR}, we describe the size of the constants $\Theta, \vartheta$ in conditional inequalities \cref{SARRH} obtained by Tsang's argument in \cite{Ts1986}.
	Assuming the Riemann Hypothesis, we can prove by Tsang's method the following inequality:
	as $T \rightarrow +\infty$, $h \rightarrow 0$,
	\begin{align}
		\label{CRKI}
		&\sup_{t \in [T, 2T]}\{ \pm(S(t + h) - S(t)) \}\\
		&\geq \l(1 - O\l( \frac{C^{k}}{\log{T}} + h C^{k} \r)\r)\frac{1}{2^{1/(2k+1)}\pi}\l( \frac{(2k)!}{k!} \r)^{1/2k}
		\l( \int_{0}^{\frac{h \log{T}}{k+1/2}}\sin^{2}\l( \frac{u}{2} \r) \frac{du}{u} \r)^{1/2}\\
		&\quad- O\l( \frac{1}{\log{T}} + h \r)
	\end{align}
	for any $k \in \ZZ_{\geq 1}$ with $k \leq c \min\{ \log{\log{T}}, \log\l( \frac{1}{h} \r) \}$.
	Below we sketch the proof of this inequality.

	First, we choose the weight function in Lemma 5 in \cite{Ts1986} as
	$V(z) = \frac{\sin z}{z}$.
	This function does not satisfy the condition of Lemma 5 in \cite{Ts1986},
	but we can easily check that this function works.

	Secondly, assuming the Riemann Hypothesis, we can show the inequality
	\begin{align}
		\sup_{u \in [T / 2, 3T]}\{\pm(S(t + h) - S(t))\}
		\geq \pm W(t) - O\l( \frac{1}{\log{T}} + h \r)
	\end{align}
	for $t \in [T, 2T]$, where
	\begin{align}
		W(t) = -\frac{2}{\pi}\Re\sum_{p^{\ell} \leq e^{\tau}, \ell \in \{ 1, 2 \}}\frac{\sin(\frac{1}{2}h\log{p^{\ell}})}{\ell p^{\ell(\frac{1}{2} + \i t)}}
	\end{align} 
	with $\tau = \frac{\log{T}}{k + 1/2} - 3\log(\frac{\log{T}}{k})$.

	Thirdly, we can prove, by using Lemma 2 in \cite{Ts1986}, that
	\begin{align}
		\int_{T}^{2T}W(t)^{2k}dt
		= \frac{T}{\pi^{2k}}\frac{(2k)!}{k!}\l( 1 + O\l( \frac{C^{k}}{\log{T}} + h C^{k} \r) \r)\l( \int_{0}^{\tau h}\sin^{2}\l( \frac{u}{2} \r)\frac{du}{u} \r)^{k},
	\end{align}
	and that
	\begin{align}
		\int_{T}^{2T}W(t)^{2k + 1}dt
		&\ll T\frac{C^{k}}{\log{T}} \l(\int_{0}^{\tau h}\sin^{2}\l(\frac{u}{2}\r) \frac{du}{u} \r)^{k + 1 / 2}.
	\end{align}
	Here, $C$ is an absolute positive constant.

	Finally, we follow the argument in the proof of Lemma 4 in \cite{Ts1986} with $R(t) = 0$ to obtain
	\begin{align}
		\max_{t \in [T, 2T]}\pm W(t)
		\geq \l(1 - O\l( \frac{C^{k}}{\log{T}} + h C^{k} \r)\r)\frac{1}{2^{1/(2k+1)}\pi}\l( \frac{(2k)!}{k!} \r)^{1/2k}
		\l( \int_{0}^{\frac{h \log{T}}{k+1/2}}\sin^{2}\l( \frac{u}{2} \r) \frac{du}{u} \r)^{\frac{1}{2}}
	\end{align}
	for any $k \in \ZZ_{\geq 1}$ with $k \leq c \min\{ \log{\log{T}}, \log\l( \frac{1}{h} \r) \}$.
	Hence, we obtain inequality \cref{CRKI}.

	By a way similar to the proof of \cref{UR_rgaps}, we can prove \cref{SARRH} with $\Theta, \vartheta$ such that
	\begin{align}
		\label{SITheta}
		\Theta
		&< \frac{1}{r^{1/2}}\frac{1}{2^{1 / (2r + 1)}\pi}\l( \frac{(2r)!}{r!} \r)^{1/2r}
		\l( \int_{0}^{\frac{4 \pi r}{2r + 1}(1+\Theta / r^{1/2})}\sin^{2}\l(\frac{u}{2}\r)\frac{du}{u} \r)^{1/2},\\
		\label{SIvtheta}
		\vartheta
		&< \frac{1}{r^{1/2}}\frac{1}{2^{1 / (2r + 1)} \pi}\l( \frac{(2r)!}{r!} \r)^{1/2r}
		\l( \int_{0}^{\frac{4 \pi r}{2r + 1}(1-\vartheta / r^{1/2})}\sin^{2}\l(\frac{u}{2}\r)\frac{du}{u} \r)^{1/2}
	\end{align}
	by using \cref{CRKI} with $k = r$ and taking the limit $T \rightarrow + \infty$.
	\Cref{TBTvt} gives approximations of $\Theta$ and $\vartheta$ satisfying \cref{SITheta,SIvtheta} for $1 \leq r \leq 20$.
	\begin{table}[h]
		\centering
		\caption{Approximations of $\Theta$, $\vartheta$ satisfying \cref{SITheta,SIvtheta}, which imply \cref{SARRH}}
		\label{TBTvt}
		\small
		\scalebox{1}[1]{
			\begin{tabular*}{140mm}{@{\extracolsep{\fill}}lcr|lcr}
				\hline
				$r$	&	$\Theta$	&	$\vartheta$	&	$r$	&	$\Theta$	&	$\vartheta$\\
				\hline
				1	&	0.394255	&	0.308149	&	11	&	0.420377	&	0.417421\\
				\hline
				2	&	0.402631	&	0.363309	&	12	&	0.4208407	&	0.418272\\
				\hline
				3	&	0.408227	&	0.385637	&	13	&	0.421236	&	0.418979\\
				\hline
				4	&	0.411824	&	0.397074	&	14	&	0.421578	&	0.419575\\
				\hline
				5	&	0.414269	&	0.403816	&	15	&	0.421876	&	0.420084\\
				\hline
				6	&	0.416023	&	0.408182	&	16	&	0.422139	&	0.420523\\
				\hline
				7	&	0.417337	&	0.411207	&	17	&	0.422372	&	0.420906\\
				\hline
				8	&	0.418355	&	0.4134103	&	18	&	0.4225802	&	0.421243\\
				\hline
				9	&	0.419167	&	0.415078	&	19	&	0.422767	&	0.4215409\\
				\hline
				10	&	0.419828	&	0.416379	&	20	&	0.422936	&	0.421806\\
				\hline
			\end{tabular*}
		}
	\end{table}
	Since $r^{-1/2}((2r)! / r!)^{1/2r}$ is increasing for $r \geq 1$,
	the right hand side of \cref{SIvtheta} is also increasing for $r \geq 1$.
	Hence, we obtain \cref{SARRH} with $\vartheta = 0.403816$ uniformly for $r \geq 5$ by calculating numerically the case $r = 5$.
	Similarly, we also obtain \cref{SARRH} with $\Theta = 0.414269$ for $r \geq 5$ by
	using the inequality $\frac{4\pi r}{2r + 1}(1 + \Theta / r^{1/2}) \geq 2\pi$.
	These are the results mentioned in \cref{Sec_CR}.

	In this note, we do not consider the unconditional explicit size of $\theta_{1}$ and $\theta_{2}$ in \cref{UR_rgaps}.
	To discuss it, we have to study implicit constants involving the zero density estimate.
	The work is related to \cite{STT2022},
	where they determined explicitly the constants $\theta_{1}, \theta_{2}$ satisfying \cref{INE_UR_rgaps} with $r = 1$.
	Moreover, they proved for such $\theta_{1}, \theta_{2}$ that inequalities \cref{INE_UR_rgaps} hold for a positive proportion of zeros.
	As a future work, the author will study the unconditional size of $\theta_{1}$ and $\theta_{2}$ for general $r$ and also consider whether we may apply the method of the present paper to consider the extension of the work \cite{STT2022} to general $r$.

	\begin{acknowledgment*}
		The author would like to thank Professor Daniel Alan Goldston for his suggestion of this problem.
		He also thanks Professor Caroline LaRoche Turnage-Butterbaugh and Professor Ade Irma Suriajaya for providing me with many valuable comments.
		The author was supported by Grant-in-Aid for JSPS Fellows (Grant Number 21J00425).
	\end{acknowledgment*}

\end{document}